\documentclass[preprint]{elsarticle}

\usepackage{amssymb,amsmath,amsthm,latexsym}
\usepackage{graphicx}
\usepackage{float}
\restylefloat{figure}
\usepackage{fullpage}

\usepackage{amsfonts,amsgen}
\usepackage{verbatim}
\usepackage{latexsym}
\usepackage{indentfirst}
\usepackage{amsmath}
\usepackage{amsfonts}

\restylefloat {figure}     
\newtheorem {theorem}{Theorem}[section]
\newtheorem {lemma}[theorem]{Lemma}

\newtheorem {definition}[theorem]{Definition}
\newtheorem {example}[theorem]{Example}

\newtheorem {conjecture}[theorem]{Conjecture}
\newtheorem {problem}[theorem]{Problem}

\input{tcilatex}

\begin{document}

\begin{frontmatter}

\title{Stability and oscillation of linear delay differential equations}

\author[label1]{John Ioannis Stavroulakis}
\author[label2]{Elena Braverman}
\address[label1]{Independent researcher}
\address[label2]{Dept. of Math. and Stats., University of
Calgary, AB, Canada T2N 1N4}

\begin{abstract}
There is a close connection between stability and oscillation of delay
differential equations. For the first-order equation
\[
x^{\prime}(t)+c(t)x(\tau(t))=0,~~t\geq 0,
\]
where $c$ is locally integrable of any sign, $\tau(t)\leq t$ is Lebesgue measurable, $\lim_{t\rightarrow\infty}\tau(t)=\infty$, we obtain sharp results, relating the speed of oscillation and stability. We thus unify the classical results of Myshkis and Lillo. We also generalise the $3/2-$stability criterion to the case of measurable parameters, improving $1+1/e$ to the sharp $3/2$ constant.
\end{abstract}


\end{frontmatter}

\section{Introduction}

The first sharp stability condition for the non-autonomous equation
with one delay term 
\begin{equation}
x^{\prime }(t)+c(t)x(\tau (t))=0,~~~t\geq 0,  \label{1}
\end{equation}
where $c:[0,+\infty )\rightarrow {{\mathbb{R}}}$ and $\tau :[0,+\infty
)\rightarrow {{\mathbb{R}}}$ are continuous, $\tau (t)\leq t$, goes back to
the work of Myshkis~\cite{28}, \cite{29}. It states that if $\func {inf}%
_{t\geq 0}c(t)>0$, $\forall t\in \lbrack 0,+\infty )$ and 
\begin{equation}
\sup _{t\in {{{{{{{\mathbb{[}}}}}}}}0,+\infty )}(t-\tau (t))~\cdot \sup
_{t\in {{{{{{{\mathbb{[}}}}}}}}0,+\infty )}c(t)<\frac{3}{2},  \label{2}
\end{equation}
all solutions of \eqref {1} tend to zero as $t\rightarrow \infty $. Lillo~%
\cite{26}, Yorke~\cite{36}, Yoneyama~\cite{33}, Gusarenko and Domoshnitskii~%
\cite{16}, successively extended criterion \eqref {2}, proving, among other
conditions, that \eqref {2} can be replaced (when $c(t)>0$, $\forall t\in
\lbrack 0,+\infty )$ and the function $t\longmapsto \int _{\tau
(t)}^{t}c(s)ds$ is continuous) by $\lim \sup _{t\rightarrow \infty }\int
_{\tau (t)}^{t}c(s)ds<\frac{3}{2}$. 

Note that the constant $\frac{3}{2}$ is sharp. According to the
results of Ladas, Sficas and Stavroulakis~\cite{25}, Malygina~\cite{27}, Gy%
\"{o}ri and Hartung~\cite{20}, it can however be improved to $\frac{\pi }{2}$
when the function $t\longmapsto \int _{\tau (t)}^{t}c(s)ds$ is sufficiently
close to a constant. 

For a non-negative coefficient $c$, any positive solution of \eqref
{1} is non-increasing and, for $\int _{0}^{\infty }c(s)ds=\infty $, tends to
zero as $t\rightarrow \infty $, and any negative solution is non-decreasing.
A function $x:[0,+\infty )\rightarrow {{\mathbb{R}}}$ is called oscillatory
if it has arbitrarily large zeros. For studying asymptotics of a solution of %
\eqref {1} with a non-negative coefficient, it is essentially sufficient to
consider oscillatory solutions.

Myshkis also studied (\cite[Theorem 12]{28}, \cite[Theorem 12]{29}),
supposing either $c(t)\geq 0$, $\forall t\in \lbrack 0,+\infty )$ or $%
c(t)\leq 0$, $\forall t\in \lbrack 0,+\infty )$ and 
\begin{equation}
\tau _{\max }:=\sup _{t\in {{{{{{{\mathbb{[}}}}}}}}0,+\infty )}(t-\tau
(t))<\infty ,  \label{19.10}
\end{equation}
solutions that have at least one zero in every interval $[t,t+\tau _{\max }]$
of length $\tau _{\max }$. Such solutions are usually referred to as rapidly
oscillating, in contrast to those for which there are always long enough
(exceeding the delay) intervals, where $x(t)$ keeps its sign. These
solutions are known as slowly oscillating. For \eqref {1} with non-negative
coefficients, existence of a slowly oscillating solution implies oscillation
of all solutions, see Agarwal et al~\cite[Theorem 2.23,p.51]{1}. As for %
\eqref {1} with non-positive coefficients, all oscillatory solutions are
rapidly oscillating, see the summary in Myshkis \cite[p.85]{28}. 

If 
\begin{equation}
\tau _{\max }\cdot ~\sup _{t\in {{{{{{{\mathbb{[}}}}}}}}0,+\infty )}|c(t)|<2,
\label{3}
\end{equation}
all rapidly oscillating solutions tend to zero as $t\rightarrow \infty $. An
example is given in Myshkis \cite{28}, showing that $2$ is the best possible
constant in \eqref {3}. 

Lillo~\cite{26}, based on previous results by Buchanan~\cite{11},
extended Myshkis' results and proved that in the case $-1\leq c(t)\leq 0$, $%
\forall t\in \lbrack 0,+\infty )$, the inequality 
\begin{equation}
\tau _{\max }<2.75+\ln 2  \label{4}
\end{equation}
implies that all oscillatory solutions of \eqref {1} tend to zero. 

We note that condition \eqref {2} of Myshkis~\cite{28}, \cite{29},
Lillo~\cite{26}, Yorke~\cite{36}, condition \eqref {3} of Myshkis~\cite{28}, 
\cite{29}, condition \eqref {4} of Lillo~\cite{26} are all three related to
certain limit-case periodic solutions of \eqref {1}. 

As for the case of \eqref {1} with a sign-changing coefficient and
multiple delays 
\begin{equation}
x^{\prime }(t)+\underset{i=1}{\overset{n}{\sum }}c_{i}(t)x(\tau _{i}(t))=0 
\label{5}
\end{equation}
the study of stability of solutions is considerably more complex. 

The first asymptotic results were obtained in the case of
\textquotedblleft a small delay\textquotedblright , by comparison with the
solutions of some simpler equations, see paragraph 10 in \cite{28}, as well
as Uvarov~\cite{32}, Driver, Sasser and Slater~\cite{13} and Driver~\cite{14}%
. Similar asymptotic results were later obtained concerning \eqref
{5} with $n=1$ by Gy\"{o}ri and Pituk \cite{18}, \cite{21}, Pituk and R\"{o}%
st~\cite{30} and with $n=2$ by Haddock and Sacker~\cite{22}, Atkinson and
Haddock~\cite{3}, Gy\"{o}ri~\cite{17}, Arino and Pituk~\cite{2}, Dibl\'{\i}k~%
\cite{12}. For $n\geq 2$, Haddock and Kuang~\cite{23} generalise the results
of \cite{3}, \cite{2}, \cite{12}, and Faria and Huang~\cite{15} generalise
those of \cite{28}, \cite{32}, \cite{13}, \cite{14}, \cite{17}. Furthermore,
Krisztin~\cite{24}, So, Yu and Chen~\cite{31}, Yoneyama and Sugie~\cite{34}, 
\cite{35}, using Yorke's method, extend the $3/2$-condition to multiple
delays. 

The technique of Azbelev, Berezansky, and Rakhmatullina \cite{4} was
applied by Gusarenko and Domoshnitskii~\cite{16}, and later by Berezansky
and Braverman \cite{6}, \cite{8}, \cite{9}. In \cite{6}, \cite{8}, \cite{9},
equation \eqref {5}  with measurable parameters was treated as a
perturbation of some canonical stable equation of type \eqref {1}, with
further application of the Bohl-Perron theorem ($W$-method). A similar
method was developed in Gy\"{o}ri, Hartung and Turi~\cite{19} and applied to
the case of a system of linear equations (the vector case with a matrix
coefficient). 

Among the various results that were obtained, is the following \cite%
{9}. Assume that $\sup _{t\in {{{{{{{\mathbb{[}}}}}}}}0,+\infty )}(t-\tau
_{i}(t))<\infty ,\quad \lim \func {inf}_{t\rightarrow \infty }\sum
_{k=1}^{n}c_{k}(t)>0$ and 
\[
\underset{t\rightarrow \infty }{\lim \sup }\int _{\underset{i=1,...,n}{\min }%
\tau _{i}(t)}^{t}\left (\underset{i=1}{\overset{n}{\sum }}c_{i}(s)\right
)ds<1+\frac{1}{e}
\]
Then all solutions of \eqref{5} tend to zero. In Berezansky and Braverman 
\cite{6}, \cite{8}, \cite{9} the extension of such results to a
sign-changing coefficient, as well as proving the $3/2$-criterion for
measurable parameters has been set as an open problem. 

In the present paper, we first prove the $3/2-$stability result in
the case of measurable parameters. Combining the techniques of Lillo \cite%
{26} and~Yorke \cite{36}, we generalise the notion of rapid oscillation and
also introduce and describe a certain real function $\Lambda $, such that,
the limit-case periodic functions corresponding to the results of Myshkis 
\cite[Theorem 12]{28}, \cite[Theorem 12]{29} and Lillo \cite{26}, can be
represented as points on its graph. We thus obtain sharp stability
conditions, relating the speed of oscillation and stability, with no sign
conditions implied on the coefficient. 

\section{The $3/2$-stability criterion for a non-negative coefficient
}

Let us introduce some relevant definitions and notations. In the
following, we assume that $c:[0,+\infty )\rightarrow {{{\mathbb{R}}}}$ is
Lebesgue measurable and locally integrable, i.e. satisfies $\int
_{0}^{t}|c(u)|du<\infty ,\forall t\in \lbrack 0,+\infty )$ and $\tau
:[0,+\infty )\rightarrow {{{\mathbb{R}}}}$ is Lebesgue measurable with 
\begin{equation}
\tau (t)\leq t\mbox {~~and~~~}\lim _{t\rightarrow \infty }\tau (t)=\infty .
\label{9a}
\end{equation}
Denote for any $t\geq 0$ 
\[
\tau _{\min }(t):=\underset{v\geq t}{\func {inf}}\tau (v)
\]
The relation \eqref {9a} guarantees that there exists $t_{1}\geq 0$ such
that $\tau _{\min }(t_{1})>0$. 

Throughout the following, we will assume 
\begin{equation}
\rho >\func {inf}\{t>0:\tau _{\min }(t)>t_{1}\}  \label{10a}
\end{equation}
In the literature, a solution of \eqref {1} is usually a function $%
x:(t_{1},t_{2}]\cup \lbrack t_{2},t_{3})\rightarrow {{{\mathbb{R}}}}$ where $%
-\infty \leq t_{1}\leq t_{2}\leq t_{3}\leq \infty $, sufficiently regular on 
$(t_{1},t_{2}]$ so that it satisfies \eqref {1} (everywhere, or almost
everywhere, depending also on the assumptions on the functions $c(t)$ and $%
\tau (t)$) on $[t_{2},t_{3})$. In this paper, we are concerned exclusively
with the behavior near $+\infty $ and not with existence/uniqueness of
solutions, so by \eqref {9a}, it suffices to study the restriction of $x$ on 
$[\tau _{\min }(t_{0}),+\infty )$ for any sufficiently large $t_{0}$. \ We
call a function $x:[\tau _{\min }(t_{0}),+\infty )\rightarrow {{{\mathbb{R}}}%
}$ , where $t_{0}\in \lbrack \rho ,+\infty )$, a solution of \eqref {1} if
it is absolutely continuous on $[\tau _{\min }(t_{0}),+\infty )$ and
satisfies \eqref {1} almost everywhere (a.e.) on $[t_{0},+\infty )$. 

First, we extend the $\frac{3}{2}$ stability criterion \eqref {2} of
Myshkis \cite{28}, \cite{29}, Lillo~\cite{26}, Yorke~\cite{36}, to the case
of measurable parameters. A generalisation of Yorke's technique~\cite[%
Proposition 4.2]{36}, \cite[Lemma 4.3]{36}, is applied. 

\begin{theorem}
\label{Theorem23.52}
Let $c(t)\geq 0,\forall t\geq 0$ and 
\begin{equation}
\sup _{t\geq \rho }\int _{\tau (t)}^{t}|c(\zeta )|d\zeta \leq \frac{3}{2}.
\label{22.44}
\end{equation}
Then all oscillatory solutions of \eqref {1} are bounded. If the inequality
in \eqref {22.44} is strict, all oscillatory solutions tend to zero. 
\end{theorem}

Define 
\[
q:=\max \left \{\underset{u\in \lbrack \rho ,+\infty )}{\sup }\int _{\tau
(u)}^{u}c(s)ds,1\right \}
\]
where $\rho $ is defined in \eqref {10a}. Without loss of generality we can
assume that at time $t$ we have 
\[
\underset{v\in \lbrack \tau _{\min }^{2}(t),t]}{\sup }|x(v)|>0,~~x(t)=0.
\]
If $\sup _{v\in \lbrack \tau _{\min }^{2}(t),t]}|x(v)|=0$ for any
sufficiently large zero $t$, we obtain for any $\varepsilon >0$, similarly
to the proof below, $\forall v\geq t,~~~|x(v)|\leq (q-\frac{1}{2}%
)\varepsilon $. That is, $x(v)=0,v\geq \tau _{\min }^{2}(t)$ and the
conclusion of the theorem holds. We note that this also follows for locally
bounded $c$, from~\cite[p.11]{4a} and \cite[Theorem B.1]{1}. 

It suffices to prove that 
\begin{equation}
\forall u\geq t,~~~|x(u)|\leq (q-\frac{1}{2})\underset{v\in \lbrack \tau
_{\min }^{2}(t),t]}{\sup }|x(v)|.  \label{22.48}
\end{equation}
If the strict inequality in \eqref {22.44} holds, $q-\frac{1}{2}\in (0,1)$,
and repeating this process ($t:=\xi _{1}$ is arbitrary) to the next zeros $%
\xi _{k}$ chosen such that $\xi _{k-1}<\tau _{\min }^{2}(\xi _{k})$, which
is possible by \eqref {9a}, one obtains 
\[
|x(u)|\leq \left (q-\frac{1}{2}\right )^{k}\underset{s\in \lbrack \tau
_{\min }^{2}(t),t]}{\sup }|x(s)|,~~u\geq \xi _{k},
\]
which would imply $\lim _{u\rightarrow \infty }x(u)=0$. 

Assume the contrary and denote, where the set in the right-hand side
is non-empty, 
\begin{align*}
z&:=\func {inf}\left \{u\geq t:|x(u)|>\left (q-\frac{1}{2}\right )\underset{%
r\in \lbrack \tau _{\min }^{2},t]}{\sup }|x(r)|\right \} \\
\tilde {t}&:=\sup \{u\in \lbrack t,z]:x(u)=0\}
\end{align*}
Obviously, $\tilde {t}\geq t$ and as $q-\frac{1}{2}\leq 1$ 
\begin{equation}
\underset{u\in \lbrack \tau _{\min }^{2}(t),z]}{\sup }|x(u)|=\underset{u\in
\lbrack \tau _{\min }^{2}(t),t]}{\sup }|x(u)|  \label{02.39}
\end{equation}
We can assume without loss of generality that $x(z)>0$ , i.e. 
\[
x(z)=\left (q-\frac{1}{2}\right )\underset{u\in \lbrack \tau _{\min
}^{2}(t),t]}{\sup }|x(u)|.
\]
By the definition of $z$, there is a decreasing sequence of $u_{n}>z$ such
that 
\[
u_{n}>u_{n+1},\quad n\in {{{{\mathbb{N}}}}},\quad \underset{n\rightarrow
\infty }{\lim }u_{n}=z,\quad x(u_{n})>x(z)>0.
\]
Then $\forall n\in N$, integrating \eqref {1}, we get 
\begin{equation}
0<x(u_{n})-x(z)=-\int _{z}^{u_{n}}c(w)x(\tau (w))dw.  \label{22.46}
\end{equation}
Since $c(w)\geq 0$ we can deduce from \eqref {22.46} that 
\[
(\forall n\in {{{\mathbb{N)}}}}~~~\exists r_{n}\in (z,u_{n}):~x(\tau
(r_{n}))<0.
\]
As the sequence $u_{n}>z$ tends to $z^{+}$, we can extract a subsequence $%
r_{k(n)}$ such that 
\[
r_{k(n)}>r_{k(n+1)}>z,~~n\in {{{{\mathbb{N}}}}},~~\underset{n\rightarrow
\infty }{\lim }r_{k(n)}=z,~~x(\tau (r_{k(n)}))<0.
\]
Also, we have 
\[
\int _{z}^{r_{k(n)}}c(w)dw>0,
\]
as $c(t)\geq 0$ and otherwise for sufficiently large $n$, we would have 
\[
\left \vert \int _{z}^{u_{n}}c(w)x(\tau (w))dw\right \vert \leq \left (%
\underset{w\in \lbrack z,u_{n}]}{ess\sup }|x(\tau (w))|\right )\int
_{z}^{u_{n}}c(w)dw=0,
\]
which contradicts \eqref {22.46}. Now, $x$ is positive on $(\tilde {t},z]$,
and by continuity, $x$ is non-negative on $[\tilde {t},r_{k(n)}]$ for
sufficiently large $n$. As$~~x(\tau (r_{k(n)}))<0$ , we have $\tau
(r_{k(n)})<\tilde {t}$. By the definition of $q$, 
\begin{equation}
\int _{\tilde {t}}^{z}c(w)dw<\int _{\tau (r_{k(n)})}^{r_{k(n)}}c(w)dw\leq q.
\label{22.47}
\end{equation}
For $v\in \lbrack \tilde {t},z]$, the following inequalities hold, (we use %
\eqref {1}, \eqref {02.39} and the definition of $q$), 
\begin{align*}
x^{\prime }(v)&\leq c(v)\min \left [|x(\tau (v))|,-x(\tau (v))\right ] \\
&\leq c(v)\min \left [\underset{u\in \lbrack \tau _{\min }^{2}(t),t]}{\sup }%
|x(u)|, 
\begin{cases}
0{{{\text { if }}}}\tau (v)\geq \tilde {t} \\[2mm] 
\int _{\tau (v)}^{\tilde {t}}c(w)|x(\tau (w))|dw{{{\text { if }}}}\tau
(v)\leq \tilde {t}%
\end{cases}
\right ] \\
&\leq c(v)\underset{u\in \lbrack \tau _{\min }^{2}(t),t]}{\sup }|x(u)|\min
\left [1, 
\begin{cases}
0{{{\text { if }}}}\tau (v)\geq \tilde {t} \\[2mm] 
\int _{\tau (v)}^{\tilde {t}}c(w)dw{{{\text { if }}}}\tau (v)\leq \tilde {t}%
\end{cases}
\right ] \\
&\leq c(v)\underset{u\in \lbrack \tau _{\min }^{2}(t),t]}{\sup }|x(u)|\min
\left [1, 
\begin{cases}
0{{{\text { if }}}}\tau (v)\geq \tilde {t} \\[2mm] 
q-\int _{\tilde {t}}^{v}c(w)dw{{{\text { if }}}}\tau (s)\leq \tilde {t}%
\end{cases}
\right ] \\
&\leq c(v)\underset{u\in \lbrack \tau _{\min }^{2}(t),t]}{\sup }|x(u)|\min
\left [1,q-\int _{\tilde {t}}^{v}c(w)dw\right ]
\end{align*}
Now we evaluate $\displaystyle x(z)=\left (q-\frac{1}{2}\right )\underset{%
u\in \lbrack \tau _{\min }^{2}(t),t]}{\sup }|x(u)|$: 
\[
x(z)=\int _{\tilde {t}}^{z}x^{\prime }(u)du\leq \underset{u\in \lbrack \tau
_{\min }^{2}(t),t]}{\sup }|x(u)|\int _{\tilde {t}}^{z}c(v)\min \left
[1,q-\int _{\tilde {t}}^{v}c(w)dw\right ]dv.
\]
Taking into account \eqref {22.47}, with the change of variable $r=\int _{%
\tilde {t}}^{v}c(w)dw$, we get 
\begin{align*}
\left (q-\frac{1}{2}\right )\underset{s\in \lbrack \tau _{\min }^{2}(t),t]}{%
\sup }|x(s)|=x(z)&\leq \underset{u\in \lbrack \tau _{\min }^{2}(t),t]}{\sup }%
|x(u)|\int _{[0,\int _{\tilde {t}}^{z}c(w)dw]}\min (1,q-r)dr \\
&<\underset{u\in \lbrack \tau _{\min }^{2}(t),t]}{\sup }|x(u)|\int
_{0}^{q}\min (1,q-r)dr \\
&=\underset{u\in \lbrack \tau _{\min }^{2}(t),t]}{\sup }|x(u)|\left (q-\frac{%
1}{2}\right ),
\end{align*}
since 
\[
\int _{0}^{q}\min [1,q-r]~dr=\int _{0}^{q-1}1dr+\int _{q-1}^{q}\left
(q-r\right )dr=\left .\int _{0}^{q-1}1dr+\left (qr-\frac{r^{2}}{2}\right
)\right \vert _{r=q-1}^{q}=q-\frac{1}{2}.
\]
The inequality $x(z)<x(z)$ is a contradiction, which completes the proof of %
\eqref {22.48}. This implies boundedness if \eqref {22.44} holds and
convergence to zero if the inequality in \eqref {22.44} is strict. 

\begin{theorem}
\label{Theorem23.53}Let $c(t)\geq 0,\forall t\geq 0$. Then all
nonoscilllatory solutions of \eqref {1} are monotone and have a limit for $%
t\rightarrow \infty $. Furthermore, if the integral $\int _{0}^{\infty
}c(v)dv$ diverges, they tend to zero as $t\rightarrow \infty $. 
\end{theorem}

Consider any nonoscillatory solution $x$ such that $x(t)>0,t\geq a.$
Assume that for the constant $b$ we have 
\begin{equation}
b>\func {inf}\{v>0:\tau _{\min }(v)\geq a\}.  \label{21.44}
\end{equation}
Now, using \eqref {1} and \eqref {21.44} 
\begin{equation}
x^{\prime }(t)=-c(t)x(\tau (t))\leq 0,t\geq b.  \label{21.47}
\end{equation}
So $x$ is nonincreasing and positive on $[b,+\infty )$. This proves that it
has a finite limit at infinity. 

Assuming 
\begin{equation}
\int _{0}^{\infty }c(v)dv=\infty  \label{21.58}
\end{equation}
we will show that $\lim _{t\rightarrow \infty }x(t)=0$. Assume the contrary,
i.e. that $\lim _{t\rightarrow \infty }x(t)>0$. Then for some $c\geq b$ ,
taking into account \eqref {9a}, we have 
\begin{equation}
x(\tau (v))\geq (1/2)\underset{t\rightarrow \infty }{\lim }x(t),v\geq c
\label{21.59}
\end{equation}
Now, using \eqref {21.47}, \eqref {21.58}, \eqref {21.59}, 
\begin{align*}
\underset{t\rightarrow \infty }{\lim }x(t)&=\underset{t\rightarrow \infty }{%
\lim }\left (x(c)+\int _{c}^{t}x^{\prime }(v)dv\right )=\underset{%
t\rightarrow \infty }{\lim }\left (x(c)-\int _{c}^{t}c(v)x(\tau (v))dv\right
) \\
&\leq \underset{t\rightarrow \infty }{\lim \func {inf}}\left (x(c)-(1/2)%
\underset{u\rightarrow \infty }{\lim }x(u)\int _{c}^{t}c(v)dv\right )=-\infty
\end{align*}
This contradicts $\lim \limits _{t\rightarrow \infty }x(t)>0$. Therefore, $%
\lim \limits _{t\rightarrow \infty }x(t)=0$ is proved.

\section{Stability in the case of sign-changing $c(t)$}

We start with extending the notion of rapid oscillation, in a way
that enables us to calculate the \textquotedblright speed\textquotedblright
\ of oscillation precisely. 

\begin{definition}
\label{Definition00.30}
We call a function $x:[t_{0},+\infty
)\rightarrow {{\mathbb{R}}},$ $\ell $-rapidly oscillating where $\ell \in
\lbrack 0,+\infty )$ if it has arbitrarily large zeros and for some $%
t_{1}\geq t_{0},\forall A,B\in \lbrack t_{1},+\infty )$ with $A<B$ the
following implication holds 
\[
\left [\forall v\in (A,B),x(v)\neq 0\right ]\Longrightarrow \int
_{A}^{B}|c(v)|dv\leq \ell .
\]
\end{definition}

In other words, $\ell $-rapid oscillation means that in any
non-oscillation segment, the integral of the absolute value of the
coefficient does not exceed a prescribed value $\ell $. 

\begin{theorem}
\label{Theorem20.54}
Let $c(t)\leq 0,t\geq 0$ and $\sup _{t\geq \rho
}\int _{\tau (t)}^{t}|c(\zeta )|d\zeta <\infty ,$where $\rho $ satisfies %
\eqref {10a}. Then all oscillatory solutions of \eqref {1} are $\sup _{t\geq
\rho }\int _{\tau (t)}^{t}|c(\zeta )|d\zeta -$rapidly oscillating.
\end{theorem}

Let $x(t)$ be a non-trivial solution of \eqref {1}. Consider two
points $\phi >\varphi $ such that $x(\varphi )=x(\phi )=0$ and $x(v)\neq
0,v\in (\varphi ,\phi )$. Without loss of generality, we can assume $%
x(v)>0,v\in (\varphi ,\phi )$. Then for $v\in (\varphi ,\phi )$ arbitrarily
close to $\phi $ 
\begin{equation}
0>x(\phi )-x(v)=-\int _{v}^{\phi }c(u)x(\tau (u))du.  \label{23.13}
\end{equation}
As $c(u)\leq 0,u\in \lbrack v,\phi ]$, by \eqref {23.13} we can assume 
\begin{equation}
\left (\forall v\in (\varphi ,\phi )\right )~~\exists \omega _{v}\in (v,\phi
):~~x(\tau (\omega _{v}))<0.  \label{20.21}
\end{equation}
Consider the sequence $v_{n}$ with $v_{0}:=\omega _{\frac{\phi +\varphi }{2}}
$ and $v_{n+1}=\omega _{\frac{v_{n}+\phi }{2}}$. Obviously, by \eqref{20.21}, 
\begin{equation}
\phi >v_{n+1}>\frac{v_{n}+\phi }{2},n\in {{\mathbb{N}}}{{\text { and}}}\lim
_{n\rightarrow \infty }v_{n}=\phi ,  \label{20.26}
\end{equation}
\begin{equation}
x(\tau (v_{n}))<0.  \label{20.23}
\end{equation}
Inequality \eqref {20.23} implies (we recall that $x(v)>0$, $v\in (\varphi
,\phi )$) 
\begin{equation}
\tau (v_{n})<\varphi .  \label{20.25}
\end{equation}
Now, using \eqref {20.25} we have 
\begin{align*}
\sup _{t\geq \rho }\int _{\tau (t)}^{t}|c(\zeta )|d\zeta &\geq \int _{\tau
(v_{n})}^{v_{n}}|c(u)|du \\
&=-\int _{v_{n}}^{\phi }|c(u)|du+\int _{\tau (v_{n})}^{\phi }|c(u)|du \\
&\geq -\int _{v_{n}}^{\phi }|c(u)|du+\int _{\varphi }^{\phi }|c(u)|du
\end{align*}
In virtue of \eqref {20.26} we can let the $v_{n}$ tend to $\phi $ and
obtain 
\[
\sup _{t\geq \rho }\int _{\tau (t)}^{t}|c(\zeta )|d\zeta \geq \int _{\varphi
}^{\phi }|c(u)|du.
\]
This concludes the proof. 

\begin{theorem}
\label{Theorem20.55}
Let $\ell \in \lbrack 0,2]$. Then all $\ell $%
-rapidly oscillating solutions of \eqref {1} are bounded. Also, if $\ell \in
\lbrack 0,2)$, all $\ell -$rapidly oscillating solutions of \eqref {1} tend
to zero. 
\end{theorem}

Let $x(t)$ be an $\ell $-rapidly oscillating solution of \eqref {1}.
If $\ell =2$ then set $q:=2$. Otherwise (if $\ell \in \lbrack 0,2)$), fix $%
q\in \left (\ell ,2\right )$. As any $\ell $-rapidly oscillating function is 
$\left (\ell +\varepsilon \right )$-rapidly oscillating, $\forall
\varepsilon >0$, we can assume that $\ell >1$ and at time $t$ we have 
\[
\sup _{v\in \lbrack \tau _{\min }^{2}(t),t]}|x(v)|>0,~~x(t)=0.
\]
If $\sup _{v\in \lbrack \tau _{\min }^{2}(t),t]}|x(v)|=0$, for any
sufficiently large zero $t$, we obtain for any $\varepsilon >0$, similarly
to the proof below, $\forall u\geq t,~~~|x(u)|\leq \max \left \{q-1,1-\frac{1%
}{2}(q-\ell )\right \}\varepsilon $. That is, $x(u)=0,u\geq \tau _{\min
}^{2}(t)$ and the conclusion of the theorem holds. We note that this also
follows for locally bounded $c$, from~\cite[p.11]{4a} and \cite[Theorem B.1]%
{1}. 

We will prove that 
\begin{equation}
|x(u)|\leq \max \left \{q-1,1-\frac{1}{2}(q-\ell )\right \}\underset{v\in
\lbrack \tau _{\min }^{2},t]}{\sup }|x(v)|,~~\forall u\geq t.  \label{00.41}
\end{equation}
Obviously if $\ell \in \lbrack 0,2)$ we have $q-1\in (0,1)$, $1-\frac{1}{2}%
(q-\ell )\in (0,1)$. Then, repeating this process ($t:=\xi _{1}$ is
arbitrary) to the next zeros $\xi _{k}$ chosen such that $\xi _{k-1}<\tau
_{\min }^{2}(\xi _{k})$, one obtains 
\[
|x(u)|\leq \left (\max \left \{q-1,1-\frac{1}{2}(q-\ell )\right \}\right
)^{k}\underset{v\in \lbrack \tau _{\min }^{2}(t),t]}{\sup }|x(v)|,~~u\geq
\xi _{k},
\]
which would imply $\lim _{u\rightarrow \infty }x(u)=0$. Assume the contrary
and define, where the set in the right-hand side is non-empty, 
\begin{align*}
z&=\func {inf}\left \{v\in \lbrack t,+\infty ):|x(v)|>\max \left (q-1,1-%
\frac{1}{2}(q-\ell )\right )\sup _{u\in \lbrack \tau _{\min
}^{2}(t),t]}|x(u)|\right \} \\
\tilde {t}&=\sup \{v\in \lbrack t,z]:x(v)=0\} \\
y&=\func {inf}\{v\in \lbrack z,+\infty ):x(v)=0\}
\end{align*}
Obviously, $\tilde {t}\geq t$ and as $\max \left (q-1,1-\frac{1}{2}(q-\ell
)\right )\leq 1$ 
\begin{equation}
\underset{u\in \lbrack \tau _{\min }^{2}(t),z]}{\sup }|x(u)|=\underset{u\in
\lbrack \tau _{\min }^{2}(t),t]}{\sup }|x(u)|  \label{02.41}
\end{equation}
As $x$ is $\ell $-rapidly oscillating, 
\begin{equation}
\int _{\tilde {t}}^{y}|c(u)|du\leq \ell \leq \frac{1}{2}(q+\ell )=q-\frac{1}{%
2}(q-\ell ).  \label{22.49}
\end{equation}
Also, for $v\in (\tilde {t},y)$ , integrating \eqref {1}, using \eqref
{02.41}, we get 
\[
|x(v)|\leq \max \left \{\underset{u\in \lbrack \tilde {t},y]}{\sup }|x(u)|,%
\underset{u\in \lbrack \tau _{\min }^{2}(t),t]}{\sup }|x(u)|\right \}\int
_{v}^{y}|c(p)|dp.
\]
Consider a point $w\in (z,y):|x(w)|=\underset{v\in \lbrack \tilde {t},y]}{%
\sup }|x(v)|$. Then, 
\begin{equation}
\max \left \{|x(w)|,\underset{u\in \lbrack \tau _{\min }^{2}(t),t]}{\sup }%
|x(u)|\right \}\int _{w}^{y}|c(p)|dp\geq |x(w)|,  \label{23.20}
\end{equation}
Taking into account that 
\[
|x(w)|>|x(z)|=\max \left (q-1,1-\frac{1}{2}(q-\ell )\right )\sup _{u\in
\lbrack \tau _{\min }^{2}(t),t]}|x(u)|\geq (q-1)\sup _{u\in \lbrack \tau
_{\min }^{2}(t),t]}|x(u)|,
\]
using \eqref {23.20} we have 
\begin{equation}
\int _{w}^{y}|c(u)|du>q-1.  \label{23.21}
\end{equation}
Using \eqref {22.49}, \eqref {23.21} one has 
\begin{equation}
\int _{\tilde {t}}^{z}|c(u)|~du\leq \int _{\tilde {t}}^{y}|c(u)|~du-\int
_{w}^{y}|c(u)|~du<1-(1/2)(q-\ell )  \label{23.22}
\end{equation}
Also, $\forall v\in \lbrack \tilde {t},z]$, we obtain 
\[
|x^{\prime }(v)|\leq |c(v)||x(\tau (v))|\leq |c(v)|\underset{u\in \lbrack
\tau _{\min }^{2}(t),t]}{\sup }|x(u)|.
\]
Now we calculate 
\[
|x(z)|=\max (q-1,1-(1/2)(q-\ell ))\sup _{u\in \lbrack \tau _{\min
}^{2}(t),t]}|x(u)|.
\]
Integrating \eqref {1} and taking into account \eqref {02.41}, we get 
\[
|x(z)|\leq \underset{u\in \lbrack \tau _{\min }^{2}(t),t]}{\sup }|x(u)|\int
_{\tilde {t}}^{z}|c(p)|dp.
\]
Applying \eqref {23.22}, 
\[
|x(z)|<\underset{u\in \lbrack \tau _{\min }^{2}(t),t]}{\sup }|x(u)|\left (1-%
\frac{1}{2}(q-\ell )\right ).
\]
Thus 
\[
|x(z)|<\max \left (q-1,1-\frac{1}{2}(q-\ell )\right )\sup _{v\in \lbrack
\tau _{\min }^{2}(t),t]}|x(v)|,
\]
which contradicts to the definition of $z$. The contradiction proves \eqref
{00.41}. 

Let us introduce a function which allows to get sharp results
relating oscillation and stability. We will show that this function
describes the limit-case periodic solutions of \eqref {1}, which mark the
border from solutions of \eqref {1} that tend to zero, to solutions\ of %
\eqref {1} that have unknown, possibly unbounded behavior. We note that
given a solution of \eqref {1} that does not tend to zero, using well-known
techniques (see Myshkis \cite[paragraph 8 ]{28}, \cite[ p.139]{28}, Lillo 
\cite[p.4]{26}), one can often construct unbounded solutions of \eqref {1}
for all greater values of the delay $\tau _{\max }$ or $\sup _{t\geq \rho
}\int _{\tau (t)}^{t}|c(\zeta )|d\zeta $. 

Define the function $\Lambda :[1,2]\longrightarrow {{\mathbb{R}}}$ 
\begin{equation}
\Lambda (s):=2+s-\sqrt{2s}-\ln \left (\sqrt{2s}-1\right ),{{\text { for }}}%
s\in \lbrack 1,2]  \label{5.09}
\end{equation}

\begin{lemma}
\label{Lemma6.23}
The function $\Lambda (s)$ defined in \eqref {5.09}
is strictly decreasing, and $\Lambda ([1,2])=[2,3-\sqrt{2}-\ln \left (\sqrt{2%
}-1\right )]$. Further, the function $\sigma (s):=$ $\Lambda (s)+s,s\in
\lbrack 1,2]$ attains a global minimum at $s=9/8.$ 
\end{lemma}

We compute $\Lambda ^{\prime }$ and $\sigma ^{\prime }$ 
\begin{align*}
\Lambda ^{\prime }(s)&=\left (2+s-\sqrt{2s}-\ln \left (\sqrt{2s}-1\right
)\right )^{\prime } \\
&=1-\frac{1}{\sqrt{2s}}-\frac{1}{\sqrt{2s}}\frac{1}{\sqrt{2s}-1} \\
&=\frac{1}{\sqrt{2s}}\left (\sqrt{2s}-1-\frac{1}{\sqrt{2s}-1}\right ) \\
&=\frac{1}{\sqrt{2s}}\cdot \frac{\left (\sqrt{2s}-1\right )^{2}-1}{\sqrt{2s}%
-1}<0,~~s\in \lbrack 1,2).
\end{align*}
Also 
\begin{align*}
\sigma ^{\prime }(s)&=2-\frac{1}{\sqrt{2s}}-\frac{1}{\sqrt{2s}}\frac{1}{%
\sqrt{2s}-1} \\
&=\frac{1}{\sqrt{2s}}\left (2\sqrt{2s}-1-\frac{1}{\sqrt{2s}-1}\right ) \\
&=\frac{1}{\sqrt{2s}}\cdot \frac{\left (2\sqrt{2s}-1\right )\left (\sqrt{2s}%
-1\right )-1}{\sqrt{2s}-1} \\
&=\frac{1}{\sqrt{2s}}\cdot \frac{4s-3\sqrt{2s}}{\sqrt{2s}-1} \\
&=\frac{2\sqrt{2s}-3}{\sqrt{2s}-1}\left \{ 
\begin{array}{ll}
>0, & s\in (9/8,2], \\ 
=0, & s=9/8, \\ 
<0, & s\in \lbrack 1,9/8).%
\end{array}
\right .
\end{align*}

\begin{lemma}
\label{Lemma6.22}
For each fixed $s\in \lbrack 1,2]$, the solution $%
\psi _{s}$ of the differential equation 
\begin{equation}
\psi _{s}^{\prime }(t)=\min \left \{1,\max \left (s-t,\psi _{s}(t)\right
)\right \},t\in \lbrack 0,\Lambda (s)-1],~~\psi _{s}(0)=0  \label{01.01}
\end{equation}
is strictly increasing and satisfies $\psi _{s}$ $(\Lambda (s)-1)=1$. 
\end{lemma}

Consider the sequence of functions 
\[
\iota _{n+1}(t):=\int _{0}^{t}\min \left \{1,\max \left (s-t,\iota
_{n}(u)\right )\right \}du,t\in \lbrack 0,\Lambda (s)-1],\iota _{0}(t)\equiv
0
\]
where $n=0,1,2,...$ We have by induction, for any fixed $t\in \lbrack
0,\Lambda (s)-1],$ $\iota _{n+1}(t)\geq \iota _{n}(t)$ . By the dominated
convergence theorem ($\Lambda (s)-1\geq \iota _{n}(t),\forall t\in \lbrack
0,\Lambda (s)-1],n=0,1,2,...$), the pointwise limit function $\iota _{\infty
}(t):=\lim _{n\rightarrow \infty }\iota _{n}(t),t\in \lbrack 0,\Lambda (s)-1]
$ solves \eqref {01.01}. In view of the Lipschitzian nature of the function $%
h\longmapsto \min \left \{1,\max \left (s-t,h\right )\right \}$ , this
solution is unique. Alternatively, one could assume existence and calculate
the solution as below, subsequently verifying that it is in fact a solution.

Evidently this solution of \eqref {01.01} satisfies $\psi _{s}(t)>0,$
and $\psi _{s}^{\prime }(t)>0$ for $t\in (0,\Lambda (s)-1]$. There are two
possible cases, $s<2$ (case \textbf{i)}) and $s=2$ (case \textbf{ii)}). 

{\textbf{i) } Case $s<2$. If $t\in \lbrack 0,s-1]$, by the definition
in \eqref {01.01} (taking into account $\psi _{s}(t)\leq \int
_{0}^{t}1=t<s-(s-1)=1$) 
\[
\psi _{s}^{\prime }(t)=1,
\]
so 
\[
\psi _{s}(s-1)=s-1.
\]
Now, as $s-1<s-(s-1)=1$ there exists an interval $(s-1,s-1+\varepsilon
),\varepsilon >0$ where $\psi _{s}^{\prime }(t)=s-t$. Also the function $%
t\longrightarrow s-t$ is strictly decreasing and $\psi _{s}$ is strictly
increasing. So either $s-t>\psi _{s}(t),t\in (s-1,\Lambda (s)-1)$ or there
is a point $c_{s}$ such that 
\begin{align*}
c_{s}&:=\sup \{t\in \lbrack s-1,\Lambda (s)-1]:u\in \lbrack
s-1,t)\Longrightarrow \psi _{s}^{\prime }(u)=s-u\}\in (s-1,\Lambda (s)-1) \\
\psi _{s}(c_{s})&=s-1+\int _{s-1}^{c_{s}}(s-t)dt=s-c_{s}
\end{align*}
The following equalities are equivalent 
\begin{align*}
s-1+\int _{s-1}^{c_{s}}(s-t)dt&=s-c_{s} \\
c_{s}-1+sc_{s}-s^{2}+s-(1/2)c_{s}^{2}+(1/2)s^{2}-s+1/2&=0 \\
c_{s}+sc_{s}-(1/2)c_{s}^{2}+-(1/2)s^{2}-1/2&=0 \\
(1/2)c_{s}^{2}-(1+s)c_{s}+(1/2)(s^{2}+1)&=0
\end{align*}
Solving the equation $(1/2)c_{s}^{2}-(1+s)c_{s}+(1/2)(s^{2}+1)=0$, we have 
\[
c_{s}=(1+s)\pm \sqrt{2s}
\]
Notice that if $c_{s}$ is well-defined, $c_{s}\geq (1+s)-\sqrt{2s}$ . Also,
by the above calculation, we have $\psi _{s}((1+s)-\sqrt{2s})=s-$ $\left
((1+s)-\sqrt{2s}\right )$ $=\sqrt{2s}-1<1$. So $c_{s}$ is in fact
well-defined and equals $(1+s)-\sqrt{2s}$. Therefore (we recall the function 
$t\longrightarrow s-t$ is strictly decreasing and $\psi _{s}$ is strictly
increasing), for $t\in \lbrack (1+s)-\sqrt{2s},\Lambda (s)-1]$, 
\[
\max \left (s-t,\psi _{s}(t)\right )=\psi _{s}(t).
\]
Defining 
\[
d_{s}:=\sup \{t\in \lbrack (1+s)-\sqrt{2s},\Lambda (s)-1]:u\in \lbrack (1+s)-%
\sqrt{2s},t)\Longrightarrow \psi _{s}(u)<1\},
\]
we have by \eqref {01.01}, 
\begin{equation}
\psi _{s}^{\prime }(t)=\psi _{s}(t),t\in \lbrack (1+s)-\sqrt{2s},d_{s}].
\label{05.58}
\end{equation}
Equation \eqref {05.58} immediately implies 
\begin{equation}
\psi _{s}(t)=\left (\sqrt{2s}-1\right )\func {exp}\left (t-(1+s)+\sqrt{2s}%
\right ).  \label{06.01}
\end{equation}
Now, by \eqref {06.01} and definition \eqref {5.09}, 
\[
d_{s}=\Lambda (s)-1
\]
and 
\[
\psi _{s}(\Lambda (s)-1)=1.
\]

\textbf{ii) } Case $s=2$. Here $\Lambda (2)=2$ and $2-t>1$, $\forall
t\in \lbrack 0,\Lambda (2)-1)$ . Hence by \eqref {01.01}, 
\begin{equation}
\psi _{s}^{\prime }(t)=1,~~t\in \lbrack 0,1).  \label{06.10}
\end{equation}
Integrating \eqref {06.10} we have $\psi _{s}(t)=t,\forall t\in \lbrack 0,1].
$ 

\begin{theorem}
\label{Theorem20.56} Assume that for some fixed $s\in \lbrack 1,2]$
the following inequalities hold 
\begin{equation}
\sup _{t\geq \rho }\int _{\tau (t)}^{t}|c(\zeta )|d\zeta \leq s,\quad \ell
\leq \Lambda (s),  \label{22.51}
\end{equation}
where $\Lambda (s)$ is denoted in \eqref {5.09} and $\rho $ satisfies \eqref
{10a}. Then all $\ell $-rapidly oscillating solutions of \eqref {1} are
bounded, and if the strict inequality $\ell <\Lambda (s)$ holds, they tend
to zero. 
\end{theorem}

Without loss of generality, we can assume that $\ell >\Lambda (s)-1$. 

If $\ell <\Lambda (s)$, let the constant $\alpha \in (1+\ell
-\Lambda (s),1)$ be a solution of 
\[
\alpha -\psi _{s}(\ell -\alpha )=0.
\]
Such a solution exists, as the function $j(\alpha ):=\alpha -\psi _{s}(\ell
-\alpha ),\alpha \in \lbrack 1+\ell -\Lambda (s),1]$ is continuous and $%
j(1+\ell -\Lambda (s))j(1)=\left (\ell -\Lambda (s)\right )\left (1-\psi
_{s}(\ell -1)\right )<0$. 

Otherwise (if $\ell =\Lambda (s)$), we set $\alpha =1$. 

Without loss of generality, we can assume at time $t$ we have 
\[
\underset{u\in \lbrack \tau _{\min }^{2}(t),t]}{\sup }|x(u)|>0,~~x(t)=0.
\]
If $\sup _{u\in \lbrack \tau _{\min }^{2}(t),t]}|x(u)|=0$ for any
sufficiently large zero $t$, we obtain for any $\varepsilon >0$, similarly
to the proof below, $\forall u\geq t,~~~|x(u)|\leq \psi _{s}(\ell -\alpha
)\varepsilon $. That is, $x(u)=0,u\geq \tau _{\min }^{2}(t)$ and the
conclusion of the theorem holds. We note that this also follows for locally
bounded $c$, from~\cite[p.11]{4a} and \cite[Theorem B.1]{1}. 

It suffices to prove that 
\begin{equation}
\forall u\geq t,~~~|x(u)|\leq \psi _{s}(\ell -\alpha )\underset{\zeta \in
\lbrack \tau _{\min }^{2}(t),t]}{\sup }|x(\zeta )|.  \label{6.26}
\end{equation}
If the strict inequality $\ell <\Lambda (s)$ holds, by Lemma {\ref{Lemma6.22}%
}, we have $\psi _{s}(\ell -\alpha )\in (0,1)$, and repeating this process ($%
t:=\xi _{1}$ is arbitrary) to the next zeros $\xi _{k}$ chosen such that $%
\xi _{k-1}<\tau _{\min }^{2}(\xi _{k})$, which is possible by \eqref {9a},
one obtains 
\[
|x(u)|\leq \left (\psi _{s}(\ell -\alpha )\right )^{k}\underset{\zeta \in
\lbrack \tau _{\min }^{2}(t),t]}{\sup }|x(\zeta )|,~~u\geq \xi _{k},
\]
which would imply $\lim _{u\rightarrow \infty }x(u)=0$. 

Assume the contrary that \eqref {6.26} is not satisfied and define 
\begin{align*}
z&=\func {inf}\left \{\zeta \in \lbrack t,+\infty ):|x(\zeta )|>\psi
_{s}(\ell -\alpha )\underset{v\in \lbrack \tau _{\min }^{2}(t),t]}{\sup }%
|x(v)|\right \}, \\
\tilde {t}&=\sup \{\zeta \in \lbrack t,z]:x(\zeta )=0\}, \\
y&=\func {inf}\{\zeta \in \lbrack z,+\infty ):x(\zeta )=0\}.
\end{align*}
We have $\tilde {t}\geq t$ and as $\psi _{s}(\ell -1)\leq 1$ 
\begin{equation}
\underset{\zeta \in \lbrack \tau _{\min }^{2}(t),t]}{\sup }|x(\zeta )|=%
\underset{\zeta \in \lbrack \tau _{\min }^{2}(t),z]}{\sup }|x(\zeta )|. 
\label{02.43}
\end{equation}

As $x$ is $\ell $-rapidly oscillating, 
\begin{equation}
\int _{\tilde {t}}^{y}|c(\zeta )|d\zeta \leq \ell .  \label{22.57}
\end{equation}
Also, for $\eta \in (\tilde {t},y),$ integrating \eqref {1}, using \eqref{02.43}, 
\begin{equation}
|x(\eta )|\leq \int _{\eta }^{y}|c(\zeta )||x(\tau (\zeta ))|d\zeta \leq
\max \left (\underset{\zeta \in \lbrack \tilde {t},y]}{\sup }|x(v)|,\underset%
{\zeta \in \lbrack \tau _{\min }^{2}(t),t]}{\sup }|x(\zeta )|\right )\int
_{\eta }^{y}|c(\zeta )|d\zeta .  \label{22.58}
\end{equation}
Consider a point $w\in (z,y)$ such that $|x(w)|=\sup _{\zeta \in \lbrack 
\tilde {t},y]}|x(\zeta )|$. Then 
\begin{equation}
x(w)=\underset{\zeta \in \lbrack \tilde {t},y]}{\sup }|x(\zeta
)|>|x(z)|=\psi _{s}(\ell -\alpha )\sup _{\zeta \in \lbrack \tau _{\min
}^{2}(t),t]}|x(\zeta )|.  \label{22.59}
\end{equation}
By \eqref {22.58}, 
\begin{equation}
\max \left (|x(w)|,\underset{\zeta \in \lbrack \tau _{\min }^{2}(t),t]}{\sup 
}|x(\zeta )|\right )\int _{w}^{y}|c(\zeta )|d\zeta \geq |x(w)|.
\label{23.25}
\end{equation}
Now, using \eqref {22.59}, \eqref {23.25} we can assume 
\begin{equation}
\int _{w}^{y}|c(\zeta )|d\zeta >\psi _{s}(\ell -\alpha )  \label{23.26}
\end{equation}
From \eqref {23.26} and \eqref {22.57}, using the definition of $\alpha $ 
\begin{equation}
\int _{\tilde {t}}^{z}|c(\zeta )|~ds\leq \int _{\tilde {t}}^{y}|c(\zeta
)|~d\zeta -\int _{w}^{y}|c(\zeta )|~d\zeta <\ell -\alpha .  \label{22.6}
\end{equation}
Also, $\forall \zeta \in \lbrack \tilde {t},z]$, integrating \eqref {1} we
have (taking into account \eqref{02.43}) 
\begin{align*}
&|x^{\prime }(\zeta )|=|c(\zeta )||x(\tau (\zeta ))|\leq |c(\zeta )|\min
\left (\underset{u\in \lbrack \tau _{\min }^{2}(t),t]}{\sup }|x(u)|,\left \{ 
\begin{array}{cc}
\int _{\tau (\zeta )}^{\tilde {t}}|c(u)||x(\tau (u))|du & {{\text {if }}}%
\tau (\zeta )<\tilde {t} \\ 
\underset{u\in \lbrack \tilde {t},\zeta ]}{\sup }|x(u)| & {{\text {if }}}%
\tau (\zeta )\geq \tilde {t}%
\end{array}
\right .\right ) \\
\leq &|c(\zeta )|\min \left (\underset{u\in \lbrack \tau _{\min }^{2}(t),t]}{%
\sup }|x(u)|,\max \left [\underset{u\in \lbrack \tilde {t},\zeta ]}{\sup }%
|x(u)|,\underset{u\in \lbrack \tau _{\min }^{2}(t),t]}{\sup }|x(u)|\int
_{\tau (\zeta )}^{\tilde {t}}|c(u)|du\right ]\right ) \\
\leq &|c(\zeta )|\min \left (\underset{u\in \lbrack \tau _{\min }^{2}(t),t]}{%
\sup }|x(u)|,\max \left [\underset{u\in \lbrack \tilde {t},\zeta ]}{\sup }%
|x(u)|,\underset{u\in \lbrack \tau _{\min }^{2}(t),t]}{\sup }|x(u)|\left
(\int _{\tau (\zeta )}^{\zeta }|c(u)|du-\int _{\tilde {t}}^{\zeta
}|c(u)|du\right )\right ]\right )
\end{align*}
Also, by \eqref {22.51}, 
\[
\int _{\tau (\zeta )}^{\zeta }|c(u)|du-\int _{\tilde {t}}^{\zeta
}|c(u)|du\leq \sup _{\zeta \geq \rho }\int _{\tau (\zeta )}^{\zeta
}|c(u)|du-\int _{\tilde {t}}^{\zeta }|c(u)|du\leq s-\int _{\tilde {t}%
}^{\zeta }|c(u)|du
\]
Thus 
\begin{equation}
|x^{\prime }(\zeta )|\leq |c(\zeta )|\min \left (\underset{u\in \lbrack \tau
_{\min }^{2}(t),t]}{\sup }|x(u)|,\max \left [\underset{u\in \lbrack \tilde {t%
},\zeta ]}{\sup }|x(u)|,\underset{u\in \lbrack \tau _{\min }^{2}(t),t]}{\sup 
}|x(u)|\left (s-\int _{\tilde {t}}^{\zeta }|c(v)|dv\right )\right ]\right )
\label{22.7}
\end{equation}
Now, either $s=2$ (case \textit{i}) or $s\in \lbrack 1,2)$ (case \textit{ii}). 

\textbf{Case \textit{i)} }We have $\psi _{s}(\ell -\alpha )=\ell -\alpha $
and integrating \eqref {22.7}, using \eqref {22.6} 
\[
|x(z)|\leq \underset{u\in \lbrack \tau _{\min }^{2}(t),t]}{\sup }|x(u)|\int
_{\tilde {t}}^{z}|c(\zeta )|d\zeta <\underset{u\in \lbrack \tau _{\min
}^{2}(t),t]}{\sup }|x(u)|\left (\ell -\alpha \right )=\underset{u\in \lbrack
\tau _{\min }^{2}(t),t]}{\sup }|x(u)|\psi _{s}(\ell -\alpha )
\]
we have a contradiction. 

\textbf{Case \textit{ii)}}  We have $\ell -\alpha \in (0,\Lambda (s)-1)$.

Below, we assume that 
\begin{equation}
(1+s)-\sqrt{2s}<\int _{\tilde {t}}^{z}|c(v)|dv<\ell -\alpha <\Lambda (s)-1
\label{21.19}
\end{equation}
The proof is similar otherwise and is based on the proof of Lemma {\ref%
{Lemma6.22}}. More precisely, as intermediate steps, we get estimates of
appropriate integrals that imply the proof for $\int _{\tilde {t}%
}^{z}|c(s)|ds\in \lbrack 0,(1+s)-\sqrt{2s}]$ or $\ell -\alpha \in (0,(1+s)-%
\sqrt{2s}]$. Alternatively, one could use a change of variables for \eqref{22.7}, 
and prove that if two functions $\alpha ,\beta $ with 
$\alpha (0)=\beta (0)=0$  solve the equality \eqref {01.01} and the
corresponding inequality, respectively, then $\alpha (t)\geq \beta (t),t\in
[0,\Lambda (s)-1]$.
Here we present a direct proof. That is,
we will show that 
\begin{equation}
\underset{u\in \lbrack \tilde {t},v]}{\sup }|x(u)|\leq \underset{u\in
\lbrack \tau _{\min }^{2}(t),t]}{\sup }|x(u)|\psi _{s} \left( \int _{\tilde {t}%
}^{v}|c(u)|du \right),~v\in \lbrack \tilde {t},z].  \label{03.02}
\end{equation}

We have 
\[
\int _{\tilde {t}}^{z}|c(\zeta )|d\zeta >(1+s)-\sqrt{2s}>s-1.
\]
We can introduce $z_{1}<z_{2}$, $z_{j}\in \lbrack \tilde {t},z)$, $j=1,2,$
such that that 
\begin{equation}
\int _{\tilde {t}}^{z_{1}}|c(\zeta )|d\zeta =s-1,~z_{2}:=\func {inf}\{r\in
\lbrack \tilde {t},z]:\int _{\tilde {t}}^{r}|c(\zeta )|d\zeta =(1+s)-\sqrt{2s%
}\}.  \label{22.8}
\end{equation}
Now, $\forall v\in \lbrack \tilde {t},z_{1}]$, integrating inequality \eqref
{22.7} 
\begin{equation}
\underset{u\in \lbrack \tilde {t},v]}{\sup }|x(u)|\leq \underset{u\in
\lbrack \tau _{\min }^{2}(t),t]}{\sup }|x(u)|\int _{\tilde {t}}^{v}|c(\zeta
)|d\zeta  \label{03.01}
\end{equation}

For $v=z_{1}$, we have 
\begin{equation}
\underset{u\in \lbrack \tilde {t},z_{1}]}{\sup }|x(u)|\leq (s-1)\sup _{u\in
\lbrack \tau _{\min }^{2}(t),t]}|x(u)|.  \label{22.09}
\end{equation}
Now, as $s-1<1$ there exists an interval $(z_{1},z_{1}+\varepsilon ),$ where 
$\varepsilon >0,\int _{z_{1}}^{z_{1}+\varepsilon }|c(s)|ds$ $>0$, where 
\begin{equation}
\sup _{u\in \lbrack \tau _{\min }^{2}(t),t]}|x(u)|\left (s-\int _{\tilde {t}%
}^{v}|c(u)|du\right )>\underset{u\in \lbrack \tilde {t},v]}{\sup }%
|x(u)|,v\in (z_{1},z_{1}+\varepsilon )  \label{21.20}
\end{equation}
Integrating \eqref {22.7} on the interval $(z_{1},z_{1}+\varepsilon ]$ ,
using \eqref {22.09}, \eqref {21.20}, 
\begin{align*}
|x(v)|&\leq |x(z_{1})|+\sup _{u\in \lbrack \tau _{\min
}^{2}(t),t]}|x(u)|\int _{z_{1}}^{v}|c(q)|\min \left (1,\max \left [s-\int _{%
\tilde {t}}^{q}|c(u)|du,\frac{\underset{u\in \lbrack \tilde {t},q]}{\sup }%
|x(u)|}{\sup _{u\in \lbrack \tau _{\min }^{2}(t),t]}|x(u)|}\right ]\right )dq
\\
&\leq \underset{u\in \lbrack \tau _{\min }^{2}(t),t]}{\sup }|x(u)|\left
(s-1\right )+\underset{u\in \lbrack \tau _{\min }^{2}(t),t]}{\sup }%
|x(u)|\int _{z_{1}}^{v}|c(q)|\left (s-\int _{\tilde {t}}^{q}|c(u)|du\right
)dq,
\end{align*}
Making a change of variable 
\[
r=\int _{\tilde {t}}^{q}|c(u)|du,
\]
we get 
\begin{equation}
|x(v)|\leq \sup _{u\in \lbrack \tau _{\min }^{2}(t),t]}|x(u)|\left [s-1+\int
_{[s-1,{\int _{\tilde {t}}^{v}}|c(s)|ds]}(s-r)dr\right ]  \label{08.32}
\end{equation}

We will show that 
\begin{equation}
\underset{u\in \lbrack \tilde {t},v]}{\sup }|x(u)|<\left (s-\int _{\tilde {t}%
}^{v}|c(u)|du\right )\underset{u\in \lbrack \tau _{\min }^{2}(t),t]}{\sup }%
|x(u)|,v\in \lbrack z_{1},z_{2})  \label{23.05}
\end{equation}
Assume the contrary, i.e. that there points $v\in \lbrack z_{1},z_{2})$ such
that 
\[
\underset{u\in \lbrack \tilde {t},v]}{\sup }|x(u)|\geq \left (s-\int _{%
\tilde {t}}^{v}|c(u)|du\right )\underset{u\in \lbrack \tau _{\min }^{2}(t),t]%
}{\sup }|x(u)|
\]
and that the following constant is well-defined 
\begin{equation}
\varpi :=\func {inf}\{v\in \lbrack z_{1},z_{2}):\underset{u\in \lbrack 
\tilde {t},v]}{\sup }|x(u)|\geq \left (s-\int _{\tilde {t}%
}^{v}|c(u)|du\right )\underset{u\in \lbrack \tau _{\min }^{2}(t),t]}{\sup }%
|x(u)|\}.  \label{23.28}
\end{equation}
Notice that by the definition of $\varpi $ in \eqref{23.28} and of $z_{2}$
in \eqref {22.8}, and also \eqref {21.20}, we have $\kappa :=\int _{\tilde {t%
}}^{\varpi }|c(u)|du\in (s-1,(1+s)-\sqrt{2s})$. However, it is obvious that $%
\forall \kappa \in (s-1,(1+s)-\sqrt{2s}),$ 
\begin{align*}
\left (s-\kappa \right )-\left [s-1\right ]-\int _{[s-1,\kappa ]}(s-r)dr&= \\
(1/2)\kappa ^{2}-(1+s)\kappa +(1/2)(s^{2}+1)&= \\
\frac{1}{2}\left (\kappa -(1+s)+\sqrt{2s}\right )\left (\kappa -(1+s)-\sqrt{%
2s}\right )&>0
\end{align*}
By the above computation and \eqref {08.32} 
\begin{equation}
\underset{u\in \lbrack \tilde {t},\varpi ]}{\sup }|x(u)|<\left (s-\int _{%
\tilde {t}}^{\varpi }|c(u)|du\right )\underset{u\in \lbrack \tau _{\min
}^{2}(t),t]}{\sup }|x(u)|{{\text {\ }}}.  \label{23.06}
\end{equation}
Equation \eqref{23.06} contradicts the definition of $\varpi $ in \eqref{23.28}. 
This proves \eqref {23.05}. By \eqref {23.05} and \eqref {08.32},
using the change of variable $r=\int _{\tilde {t}}^{s}|c(u)|du,$ 
\begin{equation}
\sup _{u\in \lbrack \tilde {t},z_{2}]}|x(u)|\leq \sup _{u\in \lbrack \tau
_{\min }^{2}(t),t]}|x(u)|\left (\left [s-1\right ]+\int _{[s-1,(1+s)-\sqrt{2s%
}]}(s-r)dr\right )  \label{01.22}
\end{equation}
We now compute the integral in \eqref {01.22}, 
\begin{align*}
&\left [s-1\right ]+\int _{[s-1,(1+s)-\sqrt{2s}]}(s-r)dr & & \\
&=s-1+s\left ((1+s)-\sqrt{2s}-\left (s-1\right )\right )-(1/2)\left [\left
((1+s)-\sqrt{2s}\right )^{2}-\left (s-1\right )^{2}\right ] & & \\
&=\sqrt{2s}-1 & &
\end{align*}
Hence, \eqref {01.22} can be written as 
\begin{equation}
\sup _{u\in \lbrack \tilde {t},z_{2}]}|x(u)|\leq \sup _{u\in \lbrack \tau
_{\min }^{2}(t),t]}|x(u)|\left (\sqrt{2s}-1\right ).  \label{23.07}
\end{equation}
Now, $\forall p\in \lbrack z_{2},z],$ (we use \eqref {23.07}, \eqref {22.7}
and the definition of $z_{2}$ in \eqref {22.8}) 
\begin{align*}
|x^{\prime }(p)|&\leq |c(p)|\max \left (\left (s-\int _{\tilde {t}%
}^{p}|c(u)|du\right )\underset{u\in \lbrack \tau _{\min }^{2}(t),t]}{\sup }%
|x(u)|,\underset{u\in \lbrack \tilde {t},p]}{\sup }|x(u)|\right ) \\
&\leq |c(p)|\max \left (\left (\sqrt{2s}-1-\int _{z_{2}}^{p}|c(u)|du\right )%
\underset{u\in \lbrack \tau _{\min }^{2}(t),t]}{\sup }|x(u)|,\underset{u\in
\lbrack \tilde {t},p]}{\sup }|x(u)|\right ) \\
&\leq |c(p)|\max \left (\left (\sqrt{2s}-1\right )\underset{u\in \lbrack
\tau _{\min }^{2}(t),t]}{\sup }|x(u)|,\underset{u\in \lbrack \tilde {t}%
,z_{2}]}{\sup }|x(u)|+\int _{z_{2}}^{p}|x^{\prime }(u)|du\right ) \\
&\leq |c(p)|\max \left (\left (\sqrt{2s}-1\right )\underset{u\in \lbrack
\tau _{\min }^{2}(t),t]}{\sup }|x(u)|,\left (\sqrt{2s}-1\right )\underset{%
u\in \lbrack \tau _{\min }^{2}(t),t]}{\sup }|x(u)|+\int
_{z_{2}}^{p}|x^{\prime }(u)|du\right ) \\
&=|c(p)|\left (\left (\sqrt{2s}-1\right )\underset{u\in \lbrack \tau _{\min
}^{2}(t),t]}{\sup }|x(u)|+\int _{z_{2}}^{p}|x^{\prime }(u)|du\right ).
\end{align*}
Therefore, the following five inequalities are valid: 
\[
|x^{\prime }(p)|-|c(p)|\int _{z_{2}}^{p}|x^{\prime }(u)|du\leq |c(p)|\left (%
\sqrt{2s}-1\right )\underset{u\in \lbrack \tau _{\min }^{2}(t),t]}{\sup }%
|x(u)|,
\]
\begin{align*}
|x^{\prime }(p)|\func {exp}\left (-\int _{0}^{p}|c(u)|du\right )-|c(p)|%
\func {exp}\left (-\int _{0}^{p}|c(u)|du\right )\int _{z_{2}}^{p}|x^{\prime
}(u)|du& \\
\leq |c(p)|\func {exp}\left (-\int _{0}^{p}|c(u)|du\right )\left (\sqrt{2s}%
-1\right )\underset{u\in \lbrack \tau _{\min }^{2}(t),t]}{\sup }|x(u)|,&
\end{align*}
\begin{equation}
\left [\func {exp}\left (-\int _{0}^{p}|c(u)|du\right )\int
_{z_{2}}^{p}|x^{\prime }(u)|du\right ]^{\prime }\leq \left [-\func {exp}%
\left (-\int _{0}^{p}|c(u)|du\right )\left (\sqrt{2s}-1\right )\underset{%
u\in \lbrack \tau _{\min }^{2}(t),t]}{\sup }|x(u)|\right ]^{\prime },
\label{01.56}
\end{equation}
Integrating \eqref {01.56}, we get 
\begin{align*}
&\func {exp}\left (-\int _{0}^{p}|c(u)|du\right )\int _{z_{2}}^{p}|x^{\prime
}(u)|du \\
\leq &\left (\func {exp}\left \{-\int _{0}^{z_{2}}|c(u)|du\right \}-\func {%
exp}\left \{-\int _{0}^{p}|c(u)|du\right \}\right )\left (\sqrt{2s}-1\right )%
\underset{u\in \lbrack \tau _{\min }^{2}(t),t]}{\sup }|x(u)|,
\end{align*}
Thus 
\begin{equation}
\int _{z_{2}}^{p}|x^{\prime }(u)|du\leq \left (\func {exp}\left \{\int
_{z_{2}}^{p}|c(u)|du\right \}-1\right )\left (\sqrt{2s}-1\right )\underset{%
u\in \lbrack \tau _{\min }^{2}(t),t]}{\sup }|x(u)|.  \label{00.37}
\end{equation}
Note that by \eqref {21.19} and the definition of $z_{2}$ in \eqref {22.8}
we have 
\begin{equation}
\int _{z_{2}}^{z}|c(u)|du<\ell -\alpha -(1+s)+\sqrt{2s}  \label{00.38}
\end{equation}
Obviously, in view of \eqref {00.37}, \eqref {08.32}, \eqref {03.01}, \eqref
{23.07}, inequality \eqref {03.02} holds. Using \eqref {00.38} and \eqref
{00.37} for $p=z$, 
\begin{equation}
\int _{z_{2}}^{z}|x^{\prime }(u)|du<\left (\sqrt{2s}-1\right )\left (\func {%
exp}\left [\ell -\alpha -(1+s)+\sqrt{2s}\right ]-1\right )\underset{u\in
\lbrack \tau _{\min }^{2}(t),t]}{\sup }|x(u)|  \label{00.39}
\end{equation}
Finally, by \eqref {23.07}, \eqref {00.39}, Lemma {\ref{Lemma6.22}, }and the
definition of $z$, 
\begin{align*}
|x(z)|&=\psi _{s}(\ell -\alpha )\underset{u\in \lbrack \tau _{\min
}^{2}(t),t]}{\sup }|x(u)| \\
&\leq \int _{z_{2}}^{z}|x^{\prime }(u)|du+\underset{u\in \lbrack \tilde {t}%
,z_{2}]}{\sup }|x(u)| \\
&<\left (\sqrt{2s}-1\right )\left (\func {exp}\left [\ell -\alpha -(1+s)+%
\sqrt{2s}\right ]-1\right )\underset{u\in \lbrack \tau _{\min }^{2}(t),t]}{%
\sup }|x(u)|+\underset{u\in \lbrack \tau _{\min }^{2}(t),t]}{\sup }%
|x(u)|\left (\sqrt{2s}-1\right ) \\
&=\left (\sqrt{2s}-1\right )\left (\func {exp}\left [\ell -\alpha -(1+s)+%
\sqrt{2s}\right ]\right )\underset{u\in \lbrack \tau _{\min }^{2}(t),t]}{%
\sup }|x(u)| \\
&=\psi _{s}(\ell -\alpha )\underset{u\in \lbrack \tau _{\min }^{2}(t),t]}{%
\sup }|x(u)|
\end{align*}
So we obtain a contradiction, which concludes the proof of \eqref {6.26} and
of the theorem. 

\begin{theorem}
\label{Theorem19.14} 
Let the assumptions of Theorem {\ref%
{Theorem20.55}} (or {\ref{Theorem20.56}}) be satisfied, and also there exist
finite $C,D>0$ such that $\forall t\geq t_{1}$, 
\begin{equation}
\left \vert t-\tau _{\min }^{2}(t)\right \vert \leq C,~~\int _{\tau _{\min
}^{2}(t)}^{t}|c(z)|dz\geq D.  \label{23.08}
\end{equation}
Then all $\ell $-rapidly-oscillating solutions of \eqref {1} tend to zero
exponentially, i.e. there exist positive constants $M>0$ and $\gamma >0$ , $%
\delta >0$ such that 
\begin{equation}
|x(t)|\leq M\sup _{u\in \lbrack t_{1},t_{1}+\delta ]}|x(u)|e^{-\gamma
(t-t_{1})}.  \label{23.09}
\end{equation}
\end{theorem}

Fix the constant $\beta :=\left \lfloor \frac{\ell }{D}\right
\rfloor +1$ (obviously $\beta >\frac{\ell }{D}$), where $\left \lfloor \cdot
\right \rfloor $ denotes the floor function. Suppose that $r-p\geq C\beta
+C+1$, where $p,r\in \lbrack t_{1},+\infty )$ are arbitrary. Then $\tau
_{\min }^{2\beta }(r)\geq r-C\beta $ and \eqref {23.08} implies that 
\begin{equation}
\int _{r-C\beta }^{r}|c(z)|dz\geq D\beta >\ell  \label{21.28}
\end{equation}
By Definition {\ref{3} and }\eqref {21.28}, $x(t)$ has at least one zero in $%
(r-C\beta ,r)$, which we will denote as $w$. 

By \eqref {23.08}, $\tau _{\min }^{2}(w)>w-C-1>p$. Therefore, using
Theorem {\ref{Theorem20.55}} (or {\ref{Theorem20.56}}), noting that $p<\tau
_{\min }^{2}(w)\leq w<r$ 
\[
\underset{u\in \lbrack r,+\infty )}{\sup }|x(u)|\leq \underset{u\in \lbrack
w,+\infty )}{\sup }|x(u)|\leq d\underset{u\in \lbrack \tau _{\min }^{2}(w),w]%
}{\sup }|x(u)|\leq d\underset{u\in \lbrack p,r]}{\sup }|x(u)|,
\]
where $d:=\max \left \{q-1,1-\frac{1}{2}(q-\ell )\right \}$, (or $\psi
_{s}(\ell -\alpha )$). Next, by induction on $\left \lfloor \frac{t-\left
(t_{1}+C\beta +C+1\right )}{C\beta +C+1}\right \rfloor $, one obtains 
\[
\underset{u\in \lbrack t,+\infty )}{\sup }|x(u)|\leq \left (\sup _{u\in
\lbrack t_{1},t_{1}+C\beta +C+1]}|x(u)|\right )d^{1+\left \lfloor \frac{%
t-t_{1}}{C\beta +C+1}-1\right \rfloor }\leq \left (\sup _{u\in \lbrack
t_{1},t_{1}+C\beta +C+1]}|x(u)|\right )d^{\left (\frac{t-t_{1}}{C\beta +C+1}%
-1\right )}.
\]
Thus \eqref {23.09} is satisfied with 
\[
|x(t)|\leq M\sup _{u\in \lbrack t_{1},t_{1}+\delta ]}|x(u)|e^{-\gamma
(t-t_{1})},~~M=d^{-1},~~\gamma =-\frac{\ln d}{C\beta +C+1},\delta =C\beta
+C+1,
\]
which concludes the proof. 

\section{Discussion, examples and open problems}

We now introduce some limit-case periodic functions, solutions of %
\eqref {1}, which are described by the $\Lambda (s)$ function, and also
illustrate the sharpness of Theorem {\ref{Theorem20.56}}. 

\begin{example}
\label{Example19.12}
Consider an arbitrary $s\in \lbrack 1,2]$. 

For $s\in \lbrack 1,2)$ consider the $\Lambda (s)$-periodic function 
\[
\begin{array}{ll}
x_{s}(t):=\psi _{s}(t) & t\in \lbrack 0,\Lambda (s)-1] \\ 
x_{s}(t)=1-(t-\Lambda (s)+1) & t\in {{{{{{{\mathbb{[}}}}}}}}\Lambda
(s)-1,\Lambda (s)] \\ 
x_{s}(t)=x_{s}(t+\Lambda (s)) & t\in {{{{{{{\mathbb{R}}}}}}}}%
\end{array}%
\]
The function $x_{s}$\ is $\Lambda (s)$-rapidly oscillating and a solution of %
\eqref {1} with 
\[
|c_{s}(t)|=1,~t\geq 0,~~\sup _{t\geq \rho }\int _{\tau
_{s}(t)}^{t}|c_{s}(\zeta )|d\zeta =s.
\]
Here $c_{s}$ is a $\Lambda (s)$-periodic function defined on $[0,\Lambda (s)]
$ as 
\[
c_{s}(t)=\left \{ 
\begin{array}{ll}
1, & t\in \lbrack 0,\Lambda (s)-1), \\ 
-1, & t\in \lbrack \Lambda (s)-1,\Lambda (s)),%
\end{array}
\right .~~~~
\]
and $\tau _{s}$ satisfies 
\[
~\tau _{s}(t)=\left \{ 
\begin{array}{ll}
-1, & t\in \lbrack 0,s-1), \\ 
t-s, & t\in \lbrack s-1,\left (s+1\right )-\sqrt{2s}), \\ 
t, & t\in \lbrack \left (s+1\right )-\sqrt{2s},\Lambda (s)-1], \\ 
\Lambda (s)-1, & t\in \lbrack \Lambda (s)-1,\Lambda (s)).%
\end{array}
\right .{{\text { and }}}\tau _{s}(t+\Lambda (s))=\tau _{s}(t)+\Lambda
(s),t\in {{\mathbb{R}}}
\]
In other words $x_{s}(t)$ satisfies almost everywhere 
\[
\begin{array}{ll}
x_{s}^{\prime }(t)=x_{s}(-1), & t\in \lbrack 0,s-1] \\ 
x_{s}^{\prime }(t)=x_{s}(t-s), & t\in \lbrack s-1,\left (s+1\right )-\sqrt{2s%
}] \\ 
x_{s}^{\prime }(t)=x_{s}(t), & t\in \lbrack \left (s+1\right )-\sqrt{2s}%
,\Lambda (s)-1] \\ 
x_{s}^{\prime }(t)=-x_{s}(\Lambda (s)-1), & t\in \lbrack \Lambda
(s)-1,\Lambda (s)] \\ 
x_{s}^{\prime }(t)=x_{s}^{\prime }(t+\Lambda (s)), & t\in {{{{{{{\mathbb{R}}}%
}}}}}%
\end{array}%
\]
For $s=2$, consider the $2$-periodic function 
\[
\begin{array}{ll}
x_{2}(t):=\psi _{2}(t) & t\in \lbrack 0,1] \\ 
x_{2}(t)=2-t & t\in {{{{{{{\mathbb{[}}}}}}}}1,2] \\ 
x_{2}(t)=x_{2}(t+2) & t\in {{{{{{{\mathbb{R}}}}}}}}%
\end{array}%
\]
The function $x_{2}$\ is $2$-rapidly oscillating and a solution of \eqref {1}
with $|c_{2}(t)|=1,t\geq 0,$ $\sup _{t\geq \rho }\int _{\tau
_{2}(t)}^{t}|c_{2}(\zeta )|d\zeta =2$. Here $c_{2}$ is a $2$-periodic
function defined on $[0,2]$ as 
\[
c_{2}(t)=\left \{ 
\begin{array}{ll}
1, & t\in \lbrack 0,1), \\ 
-1, & t\in \lbrack 1,2),%
\end{array}
\right .~~~~
\]
and $\tau _{2}$ satisfies 
\[
\tau _{2}(t)=\left \{ 
\begin{array}{ll}
-1, & t\in \lbrack 0,1), \\ 
1, & t\in \lbrack 1,2).%
\end{array}
\right .{{\text { and }}}\tau _{2}(t+2)=\tau _{2}(t)+2,t\in {{\mathbb{R}}}
\]
In other words $x_{2}(t)$ satisfies almost everywhere 
\[
\begin{array}{ll}
x_{2}^{\prime }(t)=x_{2}(-1), & t\in \lbrack 0,1] \\ 
x_{2}^{\prime }(t)=-x_{2}(1), & t\in \lbrack 1,2] \\ 
x_{2}^{\prime }(t)=x_{2}^{\prime }(t+2), & t\in {{{{{{{\mathbb{R}}}}}}}}%
\end{array}%
\]
Hence the bound $\Lambda (s)$ in Theorem {\ref{Theorem20.56}} is sharp. 
\end{example}

Next, let us explore certain $\Lambda (s)-$rapidly oscillating
solutions of \eqref {1} with a non-positive coefficient $c(t)\leq 0,\forall
t\geq 0.$ 

\begin{example}
\label{Example19.28}
Consider an arbitrary $s\in \lbrack 1,2]$. 

For $s\in [1,2)$ consider the periodic function 
\[
\begin{array}{ll}
y_{s}(t):=\psi _{s}(t) & t\in \lbrack 0,\Lambda (s)-1] \\ 
y_{s}(t)=1-(t-\Lambda (s)+1) & t\in {{{{{{{\mathbb{[}}}}}}}}\Lambda
(s)-1,\Lambda (s)] \\ 
y_{s}(t)=-y_{s}(t+\Lambda (s)) & t\in {{{{{{{\mathbb{R}}}}}}}}%
\end{array}%
\]
the function $y_{s}$ is a $\Lambda (s)$-rapidly oscillating solution of %
\eqref {1} with 
\[
c_{s}(t)=-1,\forall t\geq 0,\quad \sup _{t\geq \rho }\int _{\tau
_{s}(t)}^{t}|c_{s}(\zeta )|d\zeta =s+\Lambda (s).
\]
Here $\tau _{s}$ satisfies 
\[
\tau _{s}(t)=\left \{ 
\begin{array}{ll}
-1-\Lambda (s), & t\in \lbrack 0,s-1), \\ 
t-s-\Lambda (s), & t\in \lbrack s-1,s+1-\sqrt{2s}), \\ 
t-\Lambda (s), & t\in \lbrack s+1-\sqrt{2s},\Lambda (s)-1], \\ 
-1, & t\in \lbrack \Lambda (s)-1,\Lambda (s)),%
\end{array}
\right .{{\text { and }}}\tau _{s}(t+\Lambda (s))=\tau _{s}(t)+\Lambda
(s),t\in {{\mathbb{R}}}
\]
leading to the equation 
\[
\begin{array}{ll}
y_{s}^{\prime }(t)=y_{s}(-1-\Lambda (s)), & t\in \lbrack 0,s-1], \\ 
y_{s}^{\prime }(t)=y_{s}(t-s-\Lambda (s)), & t\in \lbrack s-1,s+1-\sqrt{2s}],
\\ 
y_{s}^{\prime }(t)=y_{s}(t-\Lambda (s)), & t\in \lbrack s+1-\sqrt{2s}%
,\Lambda (s)-1], \\ 
y_{s}^{\prime }(t)=y_{s}(-1), & t\in \lbrack \Lambda (s)-1,\Lambda (s)], \\ 
y_{s}^{\prime }(t)=-y_{s}^{\prime }(t+\Lambda (s)), & t\in {{{{{{{\mathbb{R}}%
}}}}}}.%
\end{array}%
\]
For $s=2$, consider the periodic function 
\[
\begin{array}{ll}
y_{2}(t):=\psi _{2}(t) & t\in \lbrack 0,1] \\ 
y_{2}(t)=2-t & t\in {{{{{{{\mathbb{[}}}}}}}}1,2] \\ 
y_{2}(t)=-y_{2}(t+2) & t\in {{{{{{{\mathbb{R}}}}}}}}%
\end{array}%
\]
The function $y_{2}$\ is $2$-rapidly oscillating and a solution of \eqref {1}
with $c_{2}(t)=-1,t\geq 0,$ $\sup _{t\geq \rho }\int _{\tau
_{2}(t)}^{t}|c_{2}(\zeta )|d\zeta =4$. Here $\tau _{2}$ satisfies 
\[
~~~\tau _{2}(t)=\left \{ 
\begin{array}{ll}
-3, & t\in \lbrack 0,1), \\ 
-1, & t\in \lbrack 1,2).%
\end{array}
\right .{{\text { and }}}\tau _{2}(t+2)=\tau _{2}(t)+2,t\in {{\mathbb{R}}}
\]
In other words $y_{2}(t)$ satisfies almost everywhere 
\[
\begin{array}{ll}
y_{2}^{\prime }(t)=y_{2}(-3), & t\in \lbrack 0,1] \\ 
y_{2}^{\prime }(t)=y_{2}(-1), & t\in \lbrack 1,2] \\ 
y_{2}^{\prime }(t)=-y_{2}^{\prime }(t+2), & t\in {{{{{{{\mathbb{R}}}}}}}}%
\end{array}%
\]
Recalling the function $\sigma (s)=s+\Lambda (s)$ of Lemma {\ref{Lemma6.23},
we }notice that the function\ $y_{s}$, which is a variation of $x_{s}$ in
Example {\ref{Example19.12} }( we have $|y_{s}(t)|=x_{s}(t),t\in R$ ) so
that it satisfies \eqref {1} with $0\leq -c_{s}(s)\leq 1,t\geq 0$, has $\tau
_{\max }=\sigma (s)=s+\Lambda (s)$.
\end{example}

Consider the following periodic functions, which were instrumental
in the results of Myshkis \cite{28}, \cite{29} and Lillo \cite{26}, 
\begin{align*}
f(t)&=1-t,t\in \lbrack 0,3/2] \\
f(t)&=-1/2-\int _{0}^{t-3/2}(1-u)du,t\in \lbrack 3/2,5/2] \\
f(t)&=-f(t+5/2),t\in {{{\mathbb{R}}}}
\end{align*}
\begin{align*}
g(t)&=1-t,t\in \lbrack 0,9/8] \\
g(t)&=-1/8-\int _{0}^{t-9/8}(1-u)du,t\in \lbrack 9/8,13/8] \\
g(t)&=-1/2\func {exp}\left (t-13/8\right ),t\in \lbrack 13/8,13/8+\ln 2] \\
g(t)&=-g(t+\ln 2+13/8),t\in {{{\mathbb{R}}}}
\end{align*}
The function $f$ \ (which first appeared in Myshkis \cite{28}, \cite{29})
solves \eqref {1} with and $0\leq c(t)\leq 1$, $t\geq 0$ and $\tau _{\max
}=3/2$, and the function $g$ solves \eqref {1} with $0\leq -c(s)\leq 1,t\geq
0$ and $\tau _{\max }=2.75+\ln 2$. Lillo in \cite{26} proved the following: 

Let $\Gamma $ denote the class of solutions $\gamma (t)$ of \eqref
{1} with $0\leq c(t)\leq 1$, $t\geq 0$ and $\tau _{\max }=3/2$ which satisfy 
\[
\gamma (t)=
\left\{ \begin{array}{ll} (-1)^{i}M_{i}, & t\in [t_{2i-1},t_{2i}], \\
(-1)^{i}M_{i}f(t-t_{2i}), & t\in [t_{2i},t_{2i+1}], 
\end{array} \right.
\]
where 
\begin{align*}
M_{i+1}&\geq M_{i}\geq ...\geq \underset{i\rightarrow \infty }{\lim }M_{i}>0,
\\
t_{i}&\leq t_{i+1},{{{\text { and}}}}\underset{i\rightarrow \infty }{\lim }%
t_{i}=+\infty .
\end{align*}

Suppose $z(t)$ is an oscillatory solution of \eqref {1} with $\tau
_{\max }=3/2,$and $0\leq c(t)\leq 1,t\geq 0$ which does not tend to zero.
Then for some $\gamma _{z}\in \Gamma $ which depends on $z$, we have : 
\[
(\forall \varepsilon >0)~~\exists \nu >0:t>\nu \Longrightarrow |z(t)-\gamma
_{z}(t)|<\varepsilon
\]
Also, for the case of $0\leq -c(t)\leq 1$, $t\geq 0$ and $\tau _{\max
}=2.75+\ln 2$, Buchanan~\cite{10} proved that a certain class of oscillatory
solutions of \eqref {1} are, in the above sense, asymptotic to $g$.
Concerning the connection between the function $g$ and the asymptotic
behavior of \eqref {1} with $c(t)\leq 0,t\geq 0$ we refer to the various
remarks and comparison theorems in Lillo \cite{26}. 

By Lemma {\ref{Lemma6.23}, we have that }the function $\sigma
(s)=s+\Lambda (s)$ (which is equal to the $\tau _{\max }$ of $y_{s}$ in
Example {\ref{Example19.28}}) attains its minimum at $9/8$. In other words,
Lillo's constant in \cite{26} is equal to the global minimum of the delay of
our Examples $y_{s}$ 
\[
2.75+\ln 2=\underset{s\in \lbrack 1,2]}{\min }\sigma (s)
\]
and we also have equality of the corresponding limit-case periodic solutions
of \eqref {1} 
\[
g(t+1)=y_{9/8}(t+\ln 2+13/8),t\in {{\mathbb{R}}}
\]
Furthermore, we note that condition \eqref {3} and Theorem {\ref%
{Theorem20.55}}, are related to the limit case $x_{2}(t)$ (defined in
Example {\ref{Example19.12})}, as well as to the unbounded functions in
Myshkis \cite[paragraph 8]{28}. 

We therefore propose the following conjectures.

\begin{conjecture}
\label{Conjecture10.10} 
Assume that $c(t)\leq 0$, $t\geq 0$ and $%
\sup _{t\geq \rho }\int _{\tau (t)}^{t}|c(\zeta )|d\zeta <2.75+\ln 2$, where 
$\rho $ satisfies \eqref {10a}. Then all oscillatory solutions of \eqref {1}
tend to zero. 
\end{conjecture}

\begin{conjecture}
\label{Conjecture10.09}
Fix an arbitrary $s\in \lbrack 1,2]$. Under
the assumption 
\[
\sup _{t\geq \rho }\int _{\tau (t)}^{t}|c(\zeta )|d\zeta =s,~~|c(\zeta
)|\leq 1,~~\forall \zeta \geq 0,
\]
where $\rho $ satisfies \eqref {10a}, all solutions of \eqref {1} that are $%
\Lambda (s)$-rapidly oscillating and nonnegative (and do not tend to zero),\
are asymptotic (in the sense of Lillo \cite{26} and Buchanan~\cite{10}, \cite%
{11}) to $x_{s}$. 
\end{conjecture}

\begin{conjecture}
\label{Conjecture03.42}
Under the assumption 
\[
~~|c(\zeta )|\leq 1,~~\forall \zeta \geq 0,
\]
all solutions of \eqref {1} that are $2$-rapidly oscillating and nonnegative
(and do not tend to zero),\ are asymptotic (in the sense of Lillo \cite{26}
and Buchanan~\cite{10}, \cite{11}) to $x_{2}$. 
\end{conjecture}

Moreover, the following problem remains open.

\begin{problem}
\label{Problem12.15}
Using the relationship between the definition of 
$\Lambda (s)$ and the functions $\psi _{s},x_{s}$ of Lemma {\ref{Lemma6.22}, 
}Example {\ref{Example19.12},} extend the definition of $\Lambda (s)$, as
well as Theorem {\ref{Theorem20.56}, to} $s\in (1/e,1)$. 
\end{problem}

Presumably for $s$ close to $1/e$ (a constant related to oscillation
of \eqref {1}, see Myshkis \cite[Theorem 29]{28}, \cite{29}), we can only
approximate $\Lambda (s)$. The results of Myshkis \cite[Theorem 29]{28}, 
\cite{29} and Domshlak and Aliev \cite[Theorem 8]{12a} concerning the
distance between zeros could be useful in approximating $\Lambda (s)$ for $%
s\in (1/e,1)$, as well as calculating the oscillation speed of solutions of
oscillatory equations.

Plotting $\Lambda (s)$, $s\in (1,2)$ (in red, connecting the points
when the delay is equal to one and two) and extending the plot as the
constant $2$ of Theorem~\ref{Theorem20.55}, we obtain Fig.~\ref{figure1}. 
The points $(2,2)$ and $(9/8,\ln 2+13/8)$ described by the results of
Myshkis \cite[paragraph 8]{28} and Lillo \cite[p.9]{26}, are connected by
the graph of $\Lambda (s)$. We have \textquotedblright connected the
dots\textquotedblright \ between these two well-known celebrated results.

\begin{figure}[ht]
\centering
\vspace{-45mm}
\includegraphics[scale=0.55]{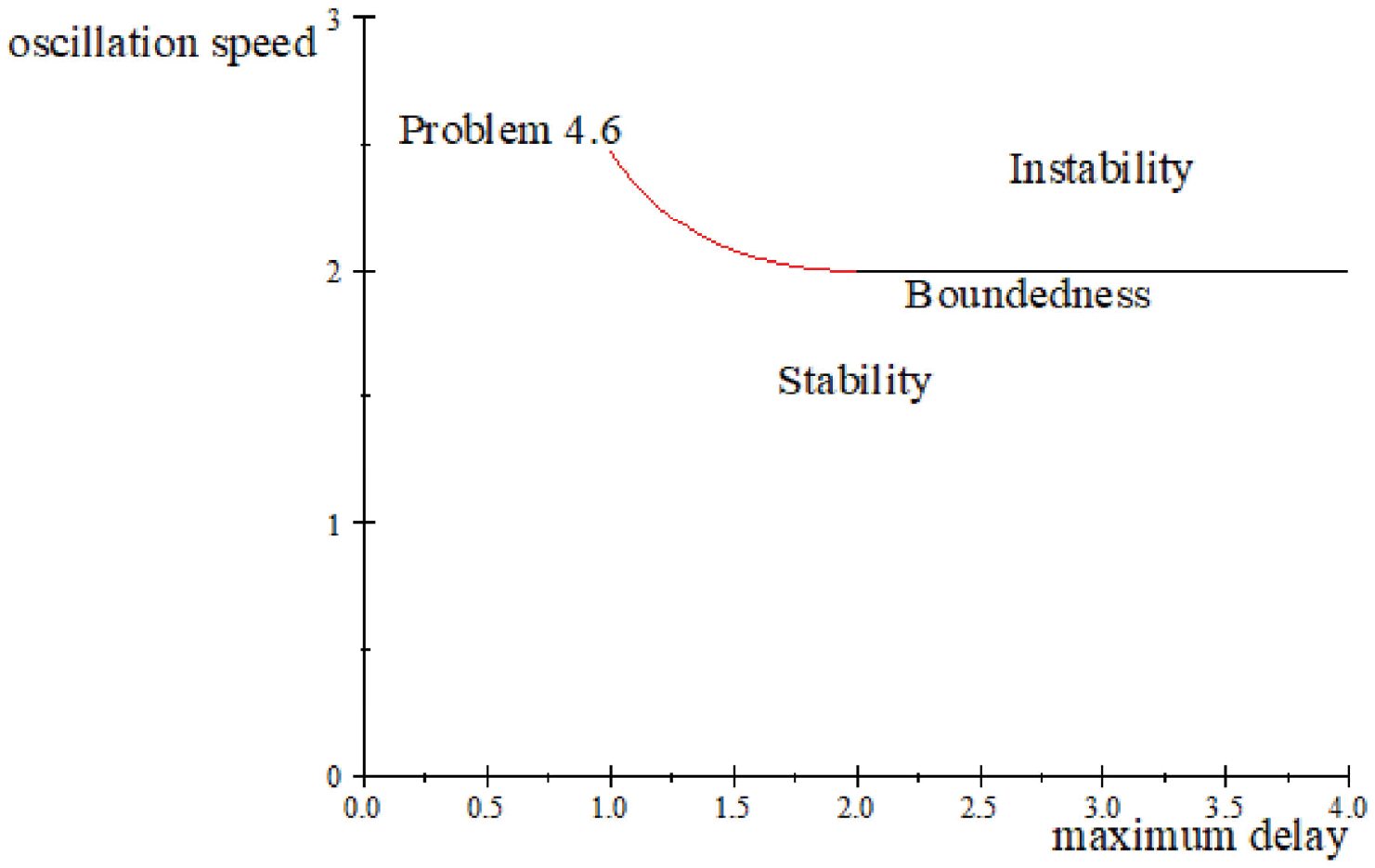} 
\caption{The graph summarizes the findings of this paper, relating the speed of
oscillation in the sense of Definition~\protect{\ref{Definition00.30}} and the maximum
delay, in the sense of  $\sup_{t\geq\rho}\int_{\tau(t)}^{t}|c(\zeta)|d\zeta$.
The periodic solutions of (\protect{\ref{1}}), $x_{s},y_{s}$ of
Examples~\protect{\ref{Example19.12}}, \protect{\ref{Example19.28}} correspond to each point of
the graph of $\Lambda(s)$. } 
\vspace{-55mm}
\label{figure1}
\end{figure}


\section*{Acknowledgment}

The second author acknowledges the support of NSERC, the grant number is \#
RGPIN-2015-05976.

\end{document}